\documentclass[reqno,10pt]{amsart}    

\usepackage{amscd,amsfonts,mathrsfs,amsthm,enumerate}
\usepackage{amssymb, amsmath}

\usepackage{color}           
\definecolor{myowncolor1}{rgb}{0,0,.75}
\definecolor{darkred}{rgb}{0.7,0,0}
\definecolor{darkgreen}{rgb}{0,0.5,0}
\definecolor{darkblue}{rgb}{0,0,.75}
\definecolor{orange}{rgb}{1,.5,0}
\definecolor{papier}{rgb}{1,1,.97}

\ifx\pdfoutput\undefined     
\usepackage{graphicx}
\else \usepackage[pdftex]{graphicx}
\usepackage[pdftex,
            colorlinks, 
            linkcolor=darkred, 
            citecolor=darkgreen, 
            urlcolor=blue]{hyperref}

\fi

\usepackage{html}            
\usepackage{ifthen}          
\usepackage{calc}            
  
\def\real     #1{{\mathbb R^{#1}}}
\def\complex  #1{{\mathbb C^{#1}}}

\def\dd       #1#2#3{{#1}_{#2#3}}

\def\dddd     #1#2#3#4#5{{#1}_{#2#3#4#5}}

\def\uu       #1#2#3{{#1}^{#2#3}}

\def\uuuu     #1#2#3#4#5{{#1}^{#2#3#4#5}}

\def\dt       {\frac{d}{dt}\,}

\def\uthetacirc #1{\theta^{^{^{~\hspace{-7pt}_\circ}}}~^{\hspace {-5pt}#1}}
\def\dthetacirc #1{\theta^{^{^{~\hspace{-7pt}_\circ}}}~_{\hspace {-7pt}#1}}

\def\dAcirc    #1{A^{^{^{~\hspace{-8pt}_\circ}}}~_{\hspace {-5pt}#1}}
\def\cM        {M^{^{^{~\hspace{-10pt}_\circ}}}}

\def\equationcolor {\color{black}}
\def\textcolor     {\color{black}}

\def\bcoleq    {\begin{equation}\equationcolor}
\def\ecoleq    {\textcolor\end{equation}}

\def\bcoleqn   {\equationcolor\begin{eqnarray}}
\def\ecoleqn   {\end{eqnarray}\textcolor}

\def\be        {$\equationcolor}
\def\ee        {\textcolor $ }

\def\bee       {$$\equationcolor}
\def\eee       {\textcolor $$}

\def\tchar     #1{\textcolor #1\equationcolor}

\def\sE        {\mathscr{E}}
\def\sF        {\mathscr{F}}

\def\picture#1#2{
\begin{latexonly}
\ifx\pdfoutput\undefined 
    \includegraphics[width=#2\hsize]
      {figures/#1.eps}
\else
    \includegraphics[width=#2\hsize]
      {figures/#1.jpg}
\fi
\end{latexonly}}

\newtheorem{theorem}{Theorem}[section]   
\newtheorem{lemma}[theorem]{Lemma}   
\newtheorem{corollary}[theorem]{Corollary}   
\newtheorem{proposition}[theorem]{Proposition}

\theoremstyle{definition}

\newcommand{\bfig}{\begin{figure}}
\newcommand{\efig}{\end{figure}}

\makeatletter
\def\pproof#1{\@ifnextchar[\opargproof
{\opargproof[\it Proof of #1.]}}
\def\opargproof[#1]{\par\noindent {\bf #1 }}

\makeatother

\begin{document}
\title
[Self-shrinkers in arbitrary codimension]
{Self-shrinkers of the mean curvature flow in arbitrary codimension}

\author[Knut Smoczyk]{\sc Knut Smoczyk}
\address{
Universit\"at Hannover,
{F}akult\"at f\"ur Mathematik und Physik,
Institut f\"ur Mathematik,    
Welfengarten 1,
30167 Hannover,
Germany}
\email{smoczyk@math.uni-hannover.de} 
\thanks{Partially supported by the DFG priority program SPP 1554, SM/78-4}

\begin{abstract}
In this paper we study self-similar solutions \be M^m\subset\real{n}\ee of the mean curvature flow
in arbitrary codimension. Self-similar curves \be\Gamma\subset\real{2}\ee have been completely 
classified by Abresch \& Langer \cite{abreschlanger} and this result can be applied to curves
\be \Gamma\subset\real{n}\ee equally well. A submanifold \be M^m\subset\real{n}\ee is called spherical,
if it is contained in a sphere. Obviously, spherical self-shrinkers of the mean curvature flow 
coincide with minimal submanifolds of the sphere. For hypersurfaces \be  M^m\subset\real{m+1},\, m\ge 2\tchar{,}\ee 
Huisken \cite{huisken1} showed that compact self-shrinkers with positive scalar mean curvature are spheres.
We will prove the following extension: A compact self-similar solution 
\be M^m\subset\real{n}\tchar{,}\, m\ge 2\tchar{,}\ee 
is spherical, if and only if the mean curvature vector \be H\ee is non-vanishing and the principal normal 
\be\nu\ee is parallel in the normal bundle. 
We also give a classification of complete noncompact self-shrinkers of that type.
\end{abstract}
\renewcommand{\subjclassname}{%
  \textup{2000} Mathematics Subject Classification}  
\subjclass[2000]{Primary 53C44; }   
\date{July 16, 2005}   
\maketitle
\section{Introduction}
In this article we study an immersion \be F:M\to\real{n}\ee of a smooth manifold \be M\ee of dimension
\be m\ee and codimension \be p=n-m\ee into euclidean space that satisfies the quasilinear elliptic system
\bcoleq\label{main 1}
H=-F^\perp,
\ecoleq
where \be H\ee denotes the mean curvature vector of the immersion and \be^\perp\ee is
the projection onto the normal bundle of \be M\tchar{.}\ee

Solutions of (\ref{main 1}) give rise to homothetically shrinking solutions \be M_t:=F(M,t)\ee of the
mean curvature flow
\bcoleqn
{F}&:&M\times[0,T)\to\real{n}\tchar{,}\nonumber\\
\frac{d}{dt}\,F(x,t)&=&H(x,t)\tchar{.}\label{main 2}
\ecoleqn
Therefore they are often called "self-shrinkers". While until recently people were mostly interested in  
mean curvature flow of hypersurfaces, this non-linear evolution equation is getting more and more
attention in higher codimension, see for instance \cite{ambson}, \cite{anciaux}, \cite{bellnov}, \cite{chenlitian}, 
\cite{tsuiwang}, \cite{smoczykwang}, \cite{ty}.

Self-similar solutions play an important role in the formation of {type$\text{-}$1} singularities
of the mean curvature flow. They are characterized by their blow-up behaviour
\bee
\sup_{M_t}|A|^2\le\frac{c}{T-t},
\eee
where \be A\ee is the second fundamental form, \be c\ee a constant and \be T\ee
denotes the blow-up time. It has been shown by Huisken \cite{huisken1} that solutions of (\ref{main 2})
forming a type-1 singularity can be homothetically rescaled so that any resulting limiting 
submanifold satisfies (\ref{main 1}).
Consequently, the classification of type-1 blow-ups is 
equivalent to the classification of self-shrinkers.

If \be M=\Gamma\subset\real{2}\ee is a curve, then all solutions of (\ref{main 1}) have been classified by Abresch and Langer
\cite{abreschlanger}. Except for the straight lines passing through the origin, 
the curvature \be k\ee is positive for all them.
In higher codimension the theorem of Abresch and Langer applies  as well.
To be precise, any self-shrinking curve \be\gamma\subset\real{n}\ee lies in a flat linear two-space \be E^2\subset\real{n}\ee and coincides
with one of the Abresch-Langer curves \be\Gamma\ee in \be E^2\tchar{,}\ee because then equation (\ref{main 1}) 
becomes an ODE of order two.

In \cite{huisken1} it was proved that a compact hypersurface \be M^{m}\subset\real{m+1},\, m\ge 2\ee
satisfying (\ref{main 1}) with positive scalar mean curvature \be H\ee is
\be S^{m}(\sqrt{m})\tchar{,}\ee 
i.e. a sphere of radius \be \sqrt{m}\tchar{.}\ee 
Later this was extended in \cite{huisken2} to noncompact hypersurfaces with \be H>0\ee 
that appear as type-1 singularities of compact hypersurfaces. There it was shown that
the blow-up limit belongs to one of the classes
\begin{itemize}
\item[i)]
\be M=\Gamma\times\real{m-1}\ee
\item[ii)]
\be M=S^{m+1-q}\bigl(\sqrt{m+1-q}\bigr)\times\real{q-1},\,q\in\{1,\dots,m\}\tchar{,}\ee
\end{itemize}
where \be\Gamma\ee are the curves found by Abresch and Langer.

If one drops the condition \be H>0\tchar{,}\ee then there are new examples, e.g. those found by Angenent
\cite{angenent}.

In higher codimension the situation becomes more complicated as the codimension increases. In particular,
even under the condition \be|H|>0\ee we get new examples. Indeed, if \be\Gamma_1,\dots,\Gamma_m\ee are
Abresch-Langer curves, then 
\bee L:=\Gamma_1\times\cdots\times\Gamma_m\subset\real{2m}\eee 
is a self-shrinker with \be |H|>0\tchar{.}\ee
Moreover, this is a Lagrangian submanifold in \be\complex{m}\ee and since the complex structure \be J\ee on
\be T\complex{m}\ee induces an isometry between the tangent and the normal bundle of \be L\tchar{,}\ee 
these examples have a flat normal bundle. Except for the case \be \Gamma_1,\dots,\Gamma_m=S^1\ee they are not 
contained in a sphere. Let us also mention that Anciaux \cite{anciaux} classified equivariant Lagrangian self-shrinkers
in \be\complex{m}\ee and that all of them satisfy \be|H|>0\ee \cite{gssz}.

On the other hand, one easily observes that spherical self-shrinkers coincide with minimal submanifolds
of the sphere and that for spherical solutions the principal normal \be\nu:=H/|H|\ee is parallel in the normal bundle.

The purpose of this article is to show that the converse holds true as well.
We state our main theorem.

\begin{theorem}\label{theorem A}
Let \be M^m\subset\real{n},\, m\ge 2\ee be a compact self-shrinker.
Then \be M\ee is spherical, if and only if \be H\neq 0\ee and  \be\nabla^\perp\nu=0\tchar{.}\ee
\end{theorem}

{\bf Remarks.}
{\it
\begin{itemize}
\item{}
It is trivial that the principal normal of a hypersurface is parallel. 
So this theorem forms a natural extension of Huisken's result mentioned above.
\item{}
In codimension two we have a torsion 1-form \be\tau\ee defined by
\bee\tau(X):=\langle\nabla_X^\perp\nu,b\rangle\tchar{,}\eee
where \be b\ee is the binormal vector \be b=J\nu\ee\ (\be J\ee denoting the complex structure on the normal bundle 
\be NM\tchar{).}\ee 
In this case 
\bee\nabla^\perp\nu=0\quad\Leftrightarrow\quad\tau=0\tchar{.}\eee
We also point out that in codimension two a parallel principal normal implies the flatness of the normal bundle, since  
the curvature \be R^\perp\ee of the normal bundle is given by \be R^\perp=d\tau\tchar{.}\ee
In higher codimensions \be\nabla^\perp\nu=0\ee is weaker than \be R^\perp=0\tchar{.}\ee
\item{}
We do not believe the condition \be\nabla^\perp\nu=0\ee is preserved under mean curvature flow.
The torsion 1-form is a generalization of the usual torsion \be\tau\ee of a biregular space curve \be\gamma\subset\real{3}\ee
and it has been shown by Altschuler \cite{altschuler} that space curves tend to develop singularities with no torsion.
Therefore, in higher codimension there is some hope that singularities will have parallel principle normal, at least under
additional geometric conditions to be discovered.
\end{itemize}
}

The idea in Huisken's original paper \cite{huisken1} was to show that the scaling invariant 
quantity \be |A|^2/|H|^2\ee is constant on a self-shrinking hypersurface. 
This will not work in our context and in fact the implications would be too strong. Instead, we look
at the quantity \be |P|^2/|H|^4\tchar{,}\ee where \be P=\langle H,A\rangle\ee is the second fundamental form w.r.t. the
mean curvature vector itself. In codimension one this coincides with \be |A|^2/|H|^2\tchar{.}\ee

One easily constructs many non-trivial examples of noncompact solutions of (\ref{main 1}). In fact, any minimal cone 
\be\Lambda\subset\real{k}\ee is
a solution and if \be M_1\subset\real{k},\,M_2\subset\real{l}\ee are solutions of (\ref{main 1}), then so is
\be M_1\times M_2\subset\real{k+l}\tchar{.}\ee In particular, a product of a minimal submanifold \be M\subset S^{n-1}\ee
with a minimal cone \be\Lambda\subset\real{k}\ee gives a noncompact solution of (\ref{main 1}) in \be\real{n+k}\ee
with parallel principal
normal \be \nu\tchar{.}\ee Since minimal cones in higher codimension need not be analytic 
(see e.g. \cite{harveylawson}), we find
non-analytic examples of noncompact self-shrinkers with parallel principal normal. 
Therefore, the general picture in the noncompact case becomes more subtle. On the other side, minimal cones never form as 
type-1 singularities of compact submanifolds, since the blow-up must have uniformly bounded curvature. In the complete 
case we give the following classification.

\begin{theorem}\label{theorem B}
Let \be M^m\subset\real{n}\ee be a complete connected self-shrinker with \be H\neq 0\ee and parallel principal normal. Suppose
further that \be M\ee has uniformly bounded geometry, i.e. there exist constants \be c_k\ee such that \be|(\nabla)^kA|\le c_k\ee
holds uniformly on \be M\ee for any \be k\ge 0\tchar{.}\ee Then \be M\ee must belong to one of the following classes.
\bee (i)\quad M^m=\Gamma\times\real{m-1}\tchar{,}\quad\quad\quad(ii)\quad M^m=\tilde M^r\times\real{m-r}\tchar{.}\eee
Here, \be\Gamma\ee is one of the curves found by Abresch and Langer and 
\be\tilde M^r\ee is a complete minimal submanifold of the
sphere \be S^{n-m+r-1}\left(\sqrt{r}\right)\subset\real{n-m+r}\tchar{,}\ee where  \be 0< r=\text{rank}(A^\nu)\le m\ee 
denotes the rank of the principal second fundamental form \be A^\nu=\langle \nu,A\rangle\tchar{.}\ee
\end{theorem}

{\bf Remarks.}
{\it
\begin{itemize}
\item{}
In codimension one, any self-shrinker
is analytic and this was exploited in \cite{huisken2} to classify the noncompact solutions with \be H\neq 0\tchar{.}\ee
In higher codimension, as pointed out above, this might not be the case. It turns out, analyticity is not needed in the proof
of the splitting theorem \ref{theorem B}. On the other hand, theorem \ref{theorem B} implies
that in codimension two any solution of (\ref{main 1}) with \be H\neq 0\ee and parallel principle normal
is analytic. 
\item{}
The uniformly bounded geometry is essential only in the sense that we need to integrate by parts w.r.t. the Gau\ss\ kernel
and this condition will prevent us from getting boundary terms at infinity. It may actually be replaced by a weaker condition which allows even a certain growth rate at infinity for \be|(\nabla)^kA|\ee and is only needed for \be k\le 2\tchar{.}\ee
In addition, completeness and bounded geometry (here, bounded second fundamental form in \be M\cap B(0,R)\ee suffices
for some large ball of radius \be R\ee centered at the origin)
will exclude products of spherical self-shrinkers with minimal cones. 
\item{}
Any blow-up (possibly noncompact)
of a type-1 singularity forming on a compact submanifold will automatically be complete with uniformly bounded geometry,
in particular theorem \ref{theorem B} may be applied to those blow-up limits.
\end{itemize}
}

The organization of the paper is as follows. After recalling the basic geometric structure equations for submanifolds
in euclidean space in section two we exploit (\ref{main 1}) to derive various elliptic equations for
curvature quantities of self-similar solutions, in particular for those with parallel principal normal. 
This will be done in section three. In section four we give a proof of theorem \ref{theorem A}, whereas in section
five we proceed to classify noncompact solutions as stated in theorem \ref{theorem B}.

This work has been initiated while I was supported by DFG as a Heisenberg fellow at the Albert Einstein Institute
in Golm and the Max Planck Institute in Leipzig and I am indebted to my hosts, Gerhard Huisken and J\"urgen Jost.
Special thanks go to Guofang Wang for discussions on this subject.   

\section{Geometry in higher codimension}\label{section 2}
Let \be F:M^m\to \real{n}\ee be a smooth immersion of a submanifold of codimension \be p=n-m\tchar{.}\ee
We let \be(x^i)_{i=1,\dots,m}\ee denote local coordinates on $M$ and we will
always use cartesian coordinates \be(y^\alpha)_{\alpha=1,\dots,n}\ee
on \be\real{n}\tchar{.}\ee Doubled greek and latin indices are summed form \be 1\ee to \be m\ee
resp. from \be 1\ee to \be n\tchar{.}\ee
In local coordinates the differential \be dF\ee of \be F\ee is given by
\bee
dF=F^\alpha_i\frac{\partial}{\partial y^\alpha}\otimes dx^i\tchar{,}
\eee
where \be F^\alpha=y^\alpha(F)\ee and
\be F^\alpha_i=\frac{\partial F^\alpha}{\partial x^i}\tchar{.}\ee
The coefficients of the induced metric \be g_{ij}\,dx^i\otimes dx^j\ee are
\bee
g_{ij}=\langle F_i,F_j\rangle=g_{\alpha\beta}F^\alpha_iF^\beta_j\tchar{,}
\eee
where \be g_{\alpha\beta}=\delta_{\alpha\beta}\ee is the euclidean metric
in cartesian coordinates. As usual, the Christoffel symbols are
\bee
\Gamma^k_{ij}=\frac 12 g^{kl}\left(\frac{\partial g_{lj}}{\partial x^i}+
\frac{\partial g_{li}}{\partial x^j}-
\frac{\partial g_{ij}}{\partial x^l}\right)\tchar{.}
\eee
The second fundmental form \be A\in\Gamma(F^{-1}T\real{n}\otimes T^*M\otimes T^*M)\ee is defined by
\bee
A:=\nabla dF=:A^\alpha_{ij}\,\frac{\partial }{\partial y^\alpha}\otimes dx^i\otimes dx^j\tchar{.}
\eee
Here and in the following all canonically induced full connections on
bundles over \be M\ee will be denoted by \be\nabla\tchar{.}\ee 
We will also use the connection on the normal bundle which will be denoted
by \be\nabla^\perp\tchar{.}\ee

It is easy to check that in cartesian coordinates on \be\real{n}\ee we have
\bee
A^\alpha_{ij}=F^\alpha_{ij}-\Gamma^k_{ij}F^\alpha_k\tchar{,}
\eee
where \be F^\alpha_{ij}=\frac{\partial^2 F^\alpha}{\partial x^i\partial x^j}\tchar{.}\ee
By definition, \be A\ee is a section in the bundle \be F^{-1}T\real{n}\otimes T^*M\otimes T^*M\ee
and it is well known that \be A\ee is normal, i.e.
\bee
A\in\Gamma\left(NM\otimes T^*M\otimes T^*M\right)\tchar{,}
\eee
where \be NM\ee denotes the normal bundle of \be M\tchar{.}\ee

This means that
\bcoleq\label{normal}
\langle F_k,\dd Aij\rangle:=g_{\alpha\beta}F^\alpha_kA^\beta_{ij}=0,\quad\forall\,i,j,k\tchar{.}
\ecoleq
The {mean curvature vector field}
\be H=H^\alpha\frac{\partial}{\partial y^\alpha}\ee is defined by
\bee H=\uu gij\dd Aij=\uu gij A^\alpha_{ij}\frac{\partial}{\partial y^\alpha}\tchar{.}\eee
It will be convenient to rise and lower indices using the metric tensors \be \dd gij,\uu gij\tchar{,}\ee
for instance
\bee
A^i_j=\uu gki\dd Akj\tchar{.}
\eee
The Riemannian curvature tensor on the tangent bundle will be denoted by \be \dddd Rijkl\tchar{,}\ee
whereas the curvature tensor of the normal bundle, considered as a 2-form with values in
\be NM\otimes NM\tchar{,}\ee has components \be R^\perp_{ij}\tchar{.}\ee 

We summarize the equations of Gau\ss, Codazzi and Ricci in the following proposition.
\begin{proposition}
{F}or an immersion \be F:M^m\to\real{n}\ee holds true
\begin{gather}
\equationcolor
\dddd Rijkl=\langle \dd Aik,\dd Ajl\rangle-\langle\dd Ail,\dd Ajk\rangle\,\tchar{,}{\textcolor\tag{Gau\ss}}\\
\equationcolor
\nabla^\perp_i\dd Ajk=\nabla^\perp_j\dd Aik\,\tchar{,}{\textcolor\tag{Codazzi}}\\
\equationcolor
R^\perp_{ij}=\dd Aik\wedge A^k_j\,\tchar{.}{\textcolor\tag{Ricci}}
\end{gather}
\end{proposition}
We define
\bee 
\dd Pij:=\langle H,\dd Aij\rangle\tchar{,}\quad \dd Qij:=\langle A_i^k,\dd Akj\rangle\tchar{,}
\quad\dddd Sijkl:=\langle \dd Aij,\dd Akl\rangle\tchar{.}
\eee
Then by Gau\ss' equation the Ricci curvature is given by
\bee
\dd Rij=\uu gkl\dddd Rikjl=\dd Pij-\dd Qij\,\tchar{.}
\eee
Moreover, we have Simons' identity
\begin{proposition}\label{ss 8}
\bcoleqn
\nabla_k^\perp\nabla_l^\perp H
&=&\Delta^\perp\dd Akl+\dddd Rkilj\uu Aij-R_k^i\dd Ail+Q_l^i\dd Aik
-\dddd Skilj\uu Aij\,\tchar{.}\nonumber
\ecoleqn
\end{proposition}
We use the Simons identity to derive an expression for \be\langle A,(\nabla^\perp)^2 H\rangle\tchar{.}\ee
In a first step we compute
\bcoleqn
2\langle A,(\nabla^\perp)^2 H\rangle
&=&\Delta|A|^2-2|\nabla^\perp A|^2+2\dddd Rkilj\uuuu Sijkl
-2\dd Rij\uu Qij\nonumber\\
&&+2|Q|^2-2\dddd Sikjl\uuuu Sijkl\,\tchar{.}\label{ss 9a}
\ecoleqn
{F}or the normal curvature tensor \be R^\perp_{ij}=\dd Aik\wedge A^k_j\ee{} one has
\bcoleqn
|R^\perp|^2
&=&|A^{\alpha i}_lA^{\beta m}_i-A^{\beta i}_lA^{\alpha m}_i|^2\nonumber\\
&=&2|Q|^2-2\dddd Sikjl\uuuu Sijkl\,\tchar{,}\nonumber
\ecoleqn
so that by Gau\ss' equation (\ref{ss 9a}) becomes
\begin{proposition}\label{ss 12}
\bcoleq
2\langle A,(\nabla^\perp)^2 H\rangle
=\Delta|A|^2-2|\nabla^\perp A|^2+2|S|^2-2\langle P,Q\rangle+2|R^\perp|^2\tchar{.}\nonumber
\ecoleq
\end{proposition}
If the immersion satisfies \be |H|>0\tchar{,}\ee we can define the {principal normal} 
\bee
\nu:=\frac{H}{|H|}\tchar{.}
\eee
The principal normal is parallel in the normal bundle, iff
\bee
\nabla_i^\perp\nu=0\tchar{.}
\eee
This is equivalent to
\bee
\quad\nabla_i^\perp H=\nabla_i|H|\nu\tchar{.}
\eee 

\section{Self-similar solutions in arbitrary codimension}\label{section 3}
Suppose \be F:M^m\to\real{n}\ee{} is a self-shrinker, i.e.
\bcoleq
H=-F^{\perp}.\nonumber
\ecoleq
We define
\bee 
\theta:=\frac{1}{2}\,d\,|F|^2
\eee
and compute
\bcoleq\label{ss 2}
\nabla_i\theta_j=\dd gij+\langle F^\perp,\dd Aij\rangle=\dd gij-\dd Pij\tchar{.}
\ecoleq
Moreover
\bcoleqn
\nabla_i^\perp F^\perp
&=&\left(\nabla_i(F-\theta^kF_k)\right)^\perp\nonumber\\
&=&\left(F_i-\nabla_i\theta^kF_k-\theta^k\dd Aik\right)^\perp\tchar{,}\nonumber
\ecoleqn
so that
\bcoleqn
\nabla_i^\perp F^\perp=-\theta^k\dd Aik\nonumber
\ecoleqn
and
\bcoleqn
\nabla_i^\perp H=\theta^k\dd Aik\tchar{.}\label{mean}
\ecoleqn
Taking another covariant derivative we derive
\bcoleqn
\nabla_i^\perp\nabla_j^\perp F^\perp
&=&-\nabla_i^\perp\left(\theta^k\dd Ajk\right)\nonumber\\
&=&-\left(\nabla_i\theta^k\dd Ajk+\theta^k\nabla_i\dd Ajk\right)^\perp\nonumber\\
&=&-\nabla_i\theta^k\dd Ajk-\theta^k\nabla_i^\perp\dd Ajk\nonumber\\
&=&-\dd Aij-\langle F^\perp,A_i^k\rangle\dd Akj-\theta^k\nabla_k^\perp\dd Aij\,\tchar{,}\nonumber
\ecoleqn
where we have used (\ref{ss 2}) and the Codazzi equation in the last step.
In particular, since \be H=-F^\perp\tchar{,}\ee we conclude
\bcoleq
\nabla^\perp_i\nabla^\perp_jH=\dd Aij-P_i^k\dd Akj+\theta^k\nabla^\perp_k\dd Aij\label{ss 4}
\ecoleq
and
\bcoleq\label{ss 6}
\Delta^\perp H-\theta^k\nabla^\perp_k H+\uu Pij\dd Aij-H=0\,\tchar{.}
\ecoleq
{F}rom this we derive
\bcoleq\label{ss 7}
\Delta|H|^2-2|\nabla^\perp H|^2-\langle F^\top,\nabla|H|^2\rangle+2|P|^2-2|H|^2=0\,\tchar{.}
\ecoleq
{F}or a self-similar solution we may exploit the Simons identity in proposition \ref{ss 8} and (\ref{ss 4}) to deduce
\bcoleq
\Delta^\perp\dd Aij-\theta^k\nabla^\perp_k\dd Aij+Q^k_i\dd Akj+Q^k_j\dd Aki+(\dddd Rikjl-\dddd Sikjl)\uu Akl
-\dd Aij=0\label{ss 9}
\ecoleq
and
\bcoleq
\Delta|A|^2-2|\nabla^\perp A|^2-\langle F^\top,\nabla|A|^2\rangle+2|S|^2+2|R^\perp|^2-2|A|^2=0
\tchar{.}\label{ss 13}
\ecoleq

We will need the following lemmata.
\begin{lemma}\label{lemma 1}
Let \be M\subset\real{n}\ee be a self-shrinker with \be |H|>0\ee and parallel
principal normal, i.e. \be\nabla^\perp\nu=0\tchar{.}\ee Then we have
\bcoleqn
\uu Pij\dd Aij&=&\frac{|P|^2}{|H|}\,\nu\,\tchar{,}\nonumber\\
\dddd Sijkl\uu Pij\uu Pkl&=&\frac{|P|^4}{|H|^2}\nonumber
\ecoleqn
and
\bcoleqn
P_i^k\dd Akj&=&P_j^k\dd Aki\,\tchar{,}\nonumber\\
\dddd Sikjl\uu Pij\uu Pkl&=&P_i^k\dd Pkj\uu Qij\,\tchar{.}\nonumber
\ecoleqn
\end{lemma}
\begin{proof}
\be\Delta^\perp H=\Delta|H|\nu,\nabla_k^\perp H=\nabla_k|H|\nu\ee and (\ref{ss 6}) imply that
\be\uu Pij\dd Aij\ee is a multiple of \be\nu\tchar{.}\ee But then 
\bee 
\uu Pij\dd Aij=\uu Pij\langle\nu,\dd Aij\rangle\nu=\frac{|P|^2}{|H|}\nu
\eee
and
\bee
\dddd Sijkl\uu Pij\uu Pkl=\langle\dd Aij\uu Pij,\dd Akl\uu Pkl\rangle=\frac{|P|^4}{|H|^2}\tchar{.}
\eee
This proves the first two equations. To prove the third we apply Ricci's equation to \be H\ee
\bee
0=\nabla_i^\perp\nabla_j^\perp H-\nabla_j^\perp\nabla_i^\perp H=\langle H,A_j^k\rangle\dd Aki-\langle H,A_i^k\rangle\dd Akj
=P_j^k\dd Aki-P_i^k\dd Akj\tchar{.}
\eee
It then follows
\bee
\dddd Sikjl\uu Pij\uu Pkl=\langle\dd Aik\uu Pij,\dd Ajl\uu Pkl\rangle
=\langle A_i^jP^i_k,\dd Ajl\uu Pkl\rangle=\dd QilP^i_k\uu Pkl
\eee
which is the last equation in the statement of the lemma.
\end{proof}
Note, that the last two equations in lemma \ref{lemma 1} hold for any submanifold with parallel principal
normal since we did not use that \be M\ee is a self-shrinker there.

\begin{lemma}\label{lemma 2}
Let \be M\subset\real{n}\ee have parallel principal normal, then
\bee
\frac{4}{|H|^4}\langle\nabla^\perp H,\nabla^\perp\dd Aij\rangle\uu Pij=
\frac{2}{|H|}\left\langle\nabla|H|,\nabla\left(\frac{|P|^2}{|H|^4}\right)\right\rangle+\frac{4|P|^2}{|H|^6}\,|\nabla|H||^2\,\tchar{.}
\eee
\end{lemma}
\begin{proof}
{F}rom \be\nabla^\perp\nu=0\ee we deduce
\bcoleqn
\langle\nabla^\perp H,\nabla^\perp\dd Aij\rangle\uu Pij
&=&\nabla^k|H|\langle\nu,\nabla^\perp_k\dd Aij\rangle\uu Pij\nonumber\\
&=&\nabla^k|H|\nabla_k\bigl(\langle\nu,\dd Aij\rangle\bigr)\uu Pij\nonumber\\
&=&\nabla^k|H|\nabla_k\left(\frac{\dd Pij}{|H|}\right)\uu Pij\nonumber\\
&=&\frac{1}{2|H|}\langle\nabla|H|,\nabla|P|^2\rangle-\frac{|P|^2}{|H|^2}|\nabla|H||^2\nonumber
\ecoleqn
and then the desired equation follows easily.
\end{proof}

A straightforward computation shows
\bcoleqn
&&\frac{2}{|H|^4}\left|\nabla_i|H|\frac{\dd Pjk}{|H|}-|H|\nabla_i\left(\frac{\dd Pjk}{|H|}\right)\right|^2
\nonumber\\
&&=\frac{2}{|H|^2}\left|\nabla\left(\frac{P}{|H|}\right)\right|^2
-2\frac{|P|^2}{|H|^6}\,|\nabla|H||^2-\frac{2}{|H|}\left\langle\nabla|H|,\nabla\left(\frac{|P|^2}{|H|^4}\right)\right\rangle\nonumber
\ecoleqn
and
\bcoleqn
&&\frac{2}{|H|^2}\left|\nabla\left(\frac{P}{|H|}\right)\right|^2\nonumber\\
&&=\frac{2|\nabla P|^2}{|H|^4}-6\frac{|P|^2}{|H|^6}|\nabla|H||^2-\frac{2}{|H|}\left\langle\nabla|H|,\nabla\left(\frac{|P|^2}{|H|^4}\right)
\right\rangle\nonumber
\ecoleqn
so that
\bcoleqn
&&\frac{2}{|H|^4}\left|\nabla_i|H|\frac{\dd Pjk}{|H|}-|H|\nabla_i\left(\frac{\dd Pjk}{|H|}\right)\right|^2
\nonumber\\
&&=\frac{2|\nabla P|^2}{|H|^4}-8\frac{|P|^2}{|H|^6}|\nabla|H||^2-\frac{4}{|H|}\left\langle\nabla|H|,\nabla\left(\frac{|P|^2}{|H|^4}\right)
\right\rangle\label{eq 1}
\ecoleqn
\begin{lemma}\label{lemma 3}
Let \be M\subset\real{n}\ee be a self-shrinker with \be|H|>0\ee and parallel principal normal \be\nu\tchar{.}\ee
Then the following equation holds
\bcoleqn
\Delta\left(\frac{|P|^2}{|H|^4}\right)
&=&\frac{2}{|H|^4}\left|\nabla_i|H|\frac{\dd Pjk}{|H|}-|H|\nabla_i\left(\frac{\dd Pjk}{|H|}\right)\right|^2\nonumber\\
&&+\left\langle F^\top,\nabla\left(\frac{|P|^2}{|H|^4}\right)\right\rangle
-\frac{2}{|H|}\left\langle\nabla|H|,\nabla\left(\frac{|P|^2}{|H|^4}\right)\right\rangle\,\tchar{.}\nonumber
\ecoleqn
\end{lemma}
\begin{proof}
{F}irst we use (\ref{ss 6}) and (\ref{ss 9}) to compute
\bcoleqn
\Delta\dd Pij
&=&\langle F^T,\nabla\dd Pij\rangle+2\langle\nabla^\perp H,\nabla^\perp\dd Aij\rangle\nonumber\\
&&-Q_i^k\dd Pkj-Q^k_j\dd Pki+2(\dddd Sikjl-\dddd Sijkl)\uu Pkl+2\dd Pij\nonumber
\ecoleqn
and then also
\bcoleqn
\Delta|P|^2
&=&2|\nabla P|^2+\langle F^\top,\nabla|P|^2\rangle+4\langle \nabla^\perp H,\nabla^\perp\dd Aij\rangle\uu Pij\nonumber\\
&&-4P_i^k\dd Pkj\uu Qij+4(\dddd Sikjl-\dddd Sijkl)\uu Pkl\uu Pij+4|P|^2\,\tchar{.}\nonumber
\ecoleqn
Since
\bcoleqn
\Delta\left(\frac{|P|^2}{|H|^4}\right)
&=&\frac{\Delta|P|^2}{|H|^4}-2\frac{|P|^2\Delta|H|^2}{|H|^6}
-\frac{8}{|H|}\left\langle\nabla|H|,\nabla\left(\frac{|P|^2}{|H|^4}\right)\right\rangle\nonumber\\
&&-8\frac{|P|^2}{|H|^6}\,|\nabla|H||^2\nonumber
\ecoleqn
this and (\ref{ss 7}) implies
\bcoleqn
\Delta\left(\frac{|P|^2}{|H|^4}\right)
&=&\left\langle F^\top,\nabla\left(\frac{|P|^2}{|H|^4}\right)\right\rangle+2\frac{|\nabla P|^2}{|H|^4}
+\frac{4}{|H|^4}\langle \nabla^\perp H,\nabla^\perp\dd Aij\rangle\uu Pij\nonumber\\
&&-4\frac{|P|^2}{|H|^6}\,|\nabla^\perp H|^2
-\frac{8}{|H|}\left\langle\nabla|H|,\nabla\left(\frac{|P|^2}{|H|^4}\right)\right\rangle\nonumber\\
&&-8\frac{|P|^2}{|H|^6}\,|\nabla|H||^2\nonumber\\
&&-\frac{4}{|H|^4}\left(P_i^k\dd Pkj\uu Qij+(\dddd Sijkl-\dddd Sikjl)\uu Pkl\uu Pij-\frac{|P|^4}{|H|^2}\right)\,\tchar{.}\nonumber
\ecoleqn
We apply lemma \ref{lemma 1} to the last term and obtain
\bcoleqn
\Delta\left(\frac{|P|^2}{|H|^4}\right)
&=&\left\langle F^\top,\nabla\left(\frac{|P|^2}{|H|^4}\right)\right\rangle+2\frac{|\nabla P|^2}{|H|^4}
+\frac{4}{|H|^4}\langle \nabla^\perp H,\nabla^\perp\dd Aij\rangle\uu Pij\nonumber\\
&&-12\frac{|P|^2}{|H|^6}\,|\nabla|H||^2
-\frac{8}{|H|}\left\langle\nabla|H|,\nabla\left(\frac{|P|^2}{|H|^4}\right)\right\rangle\,\tchar{,}\nonumber
\ecoleqn
where we have used \be\nabla^\perp H=\nabla|H|\nu+|H|\nabla^\perp\nu=\nabla|H|\nu\ee to replace \be|\nabla^\perp H|^2\ee
by \be|\nabla|H||^2\tchar{.}\ee In a next step we apply formula (\ref{eq 1}) and lemma \ref{lemma 2} to continue
\bcoleqn
\Delta\left(\frac{|P|^2}{|H|^4}\right)
&=&\frac{2}{|H|^4}\left|\nabla_i|H|\frac{\dd Pjk}{|H|}-|H|\nabla_i\left(\frac{\dd Pjk}{|H|}\right)\right|^2\nonumber\\
&&+\left\langle F^\top,\nabla\left(\frac{|P|^2}{|H|^4}\right)\right\rangle
-\frac{2}{|H|}\left\langle\nabla|H|,\nabla\left(\frac{|P|^2}{|H|^4}\right)\right\rangle\,\tchar{.}\nonumber
\ecoleqn
\end{proof}

\section{Compact self-shrinkers}
In this section we shall assume that \be M\ee is a compact self-shrinker. We are now ready to prove theorem
\ref{theorem A}.

\begin{pproof}{theorem \ref{theorem A}}
As pointed out in the introduction, it is obvious that any spherical self-shrinker is a minimal submanifold of
the sphere, that \be H\neq 0\ee  and that the principle normal is parallel. It remains to prove the converse
holds true as well.
The strong elliptic maximum principle and lemma \ref{lemma 3} imply
\bee
\frac{|P|^2}{|H|^4}=c
\eee
for some constant \be c>0\ee and
\bee
\left|\nabla_i|H|\frac{\dd Pjk}{|H|}-|H|\nabla_i\left(\frac{\dd Pjk}{|H|}\right)\right|=0\,\tchar{.}
\eee
The Codazzi equation and \be\nabla^\perp\nu=0\ee show that
\bee
\nabla_i\left(\frac{\dd Pjk}{|H|}\right)=\nabla_i\langle\nu,\dd Ajk\rangle=\langle\nu,\nabla_i^\perp\dd Ajk\rangle
=\langle\nu,\nabla_j^\perp\dd Aik\rangle=\nabla_j\left(\frac{\dd Pik}{|H|}\right)\tchar{.}
\eee
Decomposing
\be\displaystyle\nabla_i|H|\frac{\dd Pjk}{|H|}-|H|\nabla_i\left(\frac{\dd Pjk}{|H|}\right)\ee
into
\bcoleqn
&&\nabla_i|H|\frac{\dd Pjk}{|H|}-|H|\nabla_i\left(\frac{\dd Pjk}{|H|}\right)
=\frac{1}{2}\left(\nabla_i|H|\frac{\dd Pjk}{|H|}-\nabla_j|H|\frac{\dd Pik}{|H|}\right)\nonumber\\
&&+\frac{1}{2}\left(\nabla_i|H|\frac{\dd Pjk}{|H|}+\nabla_j|H|\frac{\dd Pik}{|H|}\right)
-|H|\nabla_i\left(\frac{\dd Pjk}{|H|}\right)\nonumber
\ecoleqn
we find
\bee
\left|\nabla_i|H|\frac{\dd Pjk}{|H|}-\nabla_j|H|\frac{\dd Pik}{|H|}\right|^2=0\,\tchar{.}
\eee
Consequently
\bee
|P|^2\,|\nabla|H||^2-P^i_k\uu Pkj\nabla_i|H|\nabla_j|H|=0\,\tchar{.}
\eee
If \be\nabla|H|\neq0\ee at some point \be p\in M\tchar{,}\ee then at this point there is only one nonzero eigenvalue
of the symmetric tensor \be \dd Pij\ee and the corresponding eigenvector is \be\frac{\nabla|H|}{|\nabla|H||}\,\tchar{.}\ee
At this point we have \be|P|^2=\left(\text{trace}(P)\right)^2=|H|^4\ee which implies that the constant \be c\ee from
above is \be 1\tchar{.}\ee
Since \be\nu\ee{} is parallel, equation (\ref{ss 7}) and \be|P|^2=|H|^4\ee imply
\bcoleq\label{ss 24}
\Delta|H|-\langle F^\top,\nabla|H|\rangle+|H|^3-|H|=0
\ecoleq
and as in \cite{huisken1} partial integration gives
\bee
(m-1)\int\limits_M|H|d\mu=0
\eee
which is impossible for \be m\ge 2\ee. Hence \be\nabla^\perp H=\nabla|H|\nu=0\ee{} everywhere on
\be M\ee{}. From the equation \be\Delta F=H\ee{} one deduces
\bee
\Delta|F|^2=2m+2\langle F^\perp,H\rangle=2(m-|H|^2)\tchar{,}
\eee
where we used (\ref{main 1}) in the last step. Since \be |H|^2\ee{} is constant, the maximum principle
implies 
\bcoleq\label{ss 25}
\Delta|F|^2=0=m-|H|^2\tchar{.}
\ecoleq
So \be|F|^2\ee{} is constant, \be\theta=0\ee and \be H=-F\ee. 
This proves \be M\subset S^{n-1}(\sqrt{m})\tchar{.}\ee
\end{pproof}

\section{The complete case}
In the complete case we cannot use the maximum principle equally well as in the last section. Instead we exploit integration
by parts w.r.t. the Gau\ss\ kernel \be \rho(y):=e^{-|y|^2/2}\ee on \be\real{n}\tchar{.}\ee

\begin{lemma}\label{lemma 4}
Let \be M\subset\real{n}\ee be a complete self-shrinker with \be |H|>0\tchar{,}\ee parallel principle normal \be\nu\ee and
with bounded geometry, i.e. \be|(\nabla)^kA|^2\le c_k\ee for suitable
constants \be c_k\ee and all \be k\ge 0\tchar{.}\ee Then the following integral expression holds
\bcoleqn\label{extra 11}
&&\int\limits_M\left|\nabla\left(\frac{|P|^2}{|H|^4}\right)\right|^2|H|^2\rho d\mu\nonumber\\
&&+\,2\int\limits_M\frac{|P|^2}{|H|^6}\left|\nabla_i|H|\frac{\dd Pjk}{|H|}-|H|\nabla_i\left(\frac{\dd Pjk}{|H|}\right)\right|^2\rho d\mu
=0.
\ecoleqn
\end{lemma}
\begin{proof}
The bounded geometry of \be M\ee guarantees that partial integration w.r.t. the Gau\ss\ kernel does not yield
any boundary terms at infinity. Since \be\nabla_i\rho=-\theta_i\rho\tchar{,}\ee Lemma \ref{lemma 3} and partial 
integration of \be\left(|P|^2/|H|^2\right)\Delta(|P|^2/|H|^4)\rho d\mu\ee gives (\ref{extra 11}). 
\end{proof}
\begin{corollary}
On any self-shrinker as in Lemma \ref{lemma 4} we must have
\bcoleqn
\left|\nabla_i|H|\frac{\dd Pjk}{|H|}-|H|\nabla_i\left(\frac{\dd Pjk}{|H|}\right)\right|&=&0\nonumber\tchar{,}\\
\frac{|P|^2}{|H|^4}&=&\text{const.}\nonumber
\ecoleqn
\end{corollary}
The next lemma on Abresch-Langer curves will become important in the proof of theorem \ref{theorem B} (see part (ii) below).
\begin{lemma}\label{abl}
Suppose \be\Gamma\subset\real{2}\ee is one of the self-shrinking curves found by Abresch and Langer. There exists a constant
\be c_\Gamma>0\ee such that
\bcoleq\label{al}
k\,e^{-\frac{|F|^2}{2}}=c_\Gamma
\ecoleq
holds true on all of \be\,\Gamma\tchar{.}\ee  Moreover, the critical values \be k_0\ee
of the curvature function \be k\ee satisfy
\bcoleq\label{crit}
k_{0}e^{-k_{0}^2/2}=c_\Gamma\tchar{.}
\ecoleq
If \be\,\Gamma_1,\Gamma_2\ee are two such curves with \be c_{\Gamma_1}=c_{\Gamma_2}\tchar{,}\ee
then up to an euclidean motion \be\Gamma_1=\Gamma_2\tchar{.}\ee
In particular, if \be\Gamma_1,\Gamma_2\ee are two different self-shrinking curves and 
\be k_1, k_2\ee are critical values of the curvature of \be\Gamma_1\ee resp. of \be \Gamma_2\tchar{,}\ee then \be k_1\neq k_2\tchar{.}\ee
\end{lemma}
\begin{proof}
Let \be\nu\ee be the inner unit normal. Then (\ref{main 1}) becomes \be k=-\langle F,\nu\rangle\ee and taking the gradient on both sides
gives
\bee
\nabla_i k=k\theta_i.
\eee
This implies \be\nabla(k\,e^{-\frac{|F|^2}{2}})=0\tchar{.}\ee It is clear that \be c_\Gamma\ee is positive and hence \be k>0\tchar{.}\ee
Then at a critical point \be k_0\ee we have \be\nabla_i k=\theta_i=0\ee and consequently \be k_0=- \langle F,\nu\rangle=-|F|\tchar{.}\ee
Substituting \be |F|\ee in (\ref{al}) gives (\ref{crit}). The remainder of the lemma follows from the fact that equation (\ref{al})
is a second order ODE.
\end{proof}

\begin{pproof}{theorem \ref{theorem B}}
Proceeding as in the proof of theorem \ref{theorem A} we first conclude that either \be\nabla^\perp H=0\ee everywhere,
or \be P\ee admits only one non-zero eigenvalue which is \be|H|^2\tchar{.}\ee Let us treat both cases separately.
\begin{itemize}
\item[(i)]\be\nabla^\perp H=0\tchar{:}\ee\\
There exist several classification results for submanifolds in euclidean space with parallel mean curvature. In particular, under
surprisingly similar assumptions on the relation between the mean curvature vector \be H\ee  and the position vector \be F\tchar{,}\ee
Yano \cite{yano} classified compact submanifolds with parallel mean curvature. It seems, an analogous classification to that of Yano 
for complete submanifolds with parallel mean curvature does not exist. However, in our special situation this can be established as 
follows. First observe, that in case \be\nabla^\perp H=0\tchar{,}\ee equation (\ref{mean}) implies
\bcoleq\label{extra 9}
\theta^i\dd Aij=0\tchar{.}
\ecoleq
Moreover, from \be\nabla^\perp H=0\ee we deduce \be\nabla\dd Pij=0\ee and then with equation (\ref{ss 4})
\bcoleq\label{extra 0}
P_i^j=P_i^kP_k^j\tchar{.}
\ecoleq
Hence, if \be \lambda\ee is an eigenvalue of \be P_i^j\tchar{,}\ee then either \be\lambda=0\ee or
\be\lambda=1\tchar{.}\ee 
Let us define
\bee
(P*A)_{ij}:=P_i^k\dd Akj\tchar{.}
\eee
By Lemma \ref{lemma 1}, \be P*A\ee is symmetric.
We claim
\bcoleq\label{extra 6}
\theta^k\nabla^\perp_k\dd Aij=0\quad\tchar{\text{ and }}\dd Aij=P_i^l\dd Alj\tchar{.}
\ecoleq
To prove (\ref{extra 6}), we exploit (\ref{ss 4}), (\ref{extra 0}) and \be\nabla P=0\ee to get in a first step
\bcoleqn
\theta^k\nabla^\perp_k \left(P_i^l\dd Alj\right)
&=&\theta^k P_i^l\nabla_k^\perp\dd Alj\nonumber\\
&\overset{(\ref{ss 4})}=&P_i^lP_l^m\dd Amj-P_i^l\dd Alj\nonumber\\
&\overset{(\ref{extra 0})}=&P_i^m\dd Amj-P_i^l\dd Alj\nonumber\\
&=&0\tchar{.}\label{extra 7}
\ecoleqn
So it suffices to show \be|A|^2-|P*A|^2=0\tchar{.}\ee From (\ref{ss 4}), (\ref{extra 0}) and (\ref{extra 7}) we conclude
\bcoleqn
\theta^k\nabla_k\left(|A|^2-|P*A|^2\right)
&=&2\theta^k\langle\uu Aij,\nabla^\perp_k \dd Aij\rangle\nonumber\\
&\overset{(\ref{ss 4})}=&2\langle\uu Aij,P_i^k\dd Akj-\dd Aij\rangle\nonumber\\
&\overset{(\ref{extra 0})}=&2\langle\uu Aij,P_i^lP_l^k\dd Akj-\dd Aij\rangle\nonumber\\
&=&-2\left(|A|^2-|P*A|^2\right)\tchar{.}\label{extra 8}
\ecoleqn
If \be\theta\ee vanishes at a point \be p\in M\tchar{,}\ee then (\ref{extra 8}) implies \be\dd Aij=P_i^l\dd Alj\ee and 
\be|A|^2=|P*A|^2\ee at \be p\tchar{.}\ee 
So suppose \be\theta(p)\neq 0\tchar{.}\ee Let \be\gamma\ee be the integral curve of \be\theta\ee with
\be\gamma(0)=p\tchar{,}\ee i.e. \be d\gamma=\theta\tchar{.}\ee Then along \be\gamma\ee we define the function
\bee
f(t):=|\theta|^2(\gamma(t))
\eee
and obtain
\bee
\dt f=d\gamma(\nabla|\theta|^2)=\theta^k\nabla_k|\theta|^2=2\theta^k\theta^l\nabla_k\theta_l\tchar{.}
\eee
{F}rom (\ref{ss 2}), \be F^\perp=-H\ee and  (\ref{extra 9}) we see
\bee
\nabla_i\theta_j=\dd gij-\dd Pij,\quad\theta^i\nabla_i\theta_j=\theta_j\tchar{,}
\eee
so that
\bcoleq\label{extra 10}
\dt f=2f\tchar{.}
\ecoleq
Consequently
\bee
|\theta|^2(\gamma(t))=|\theta|^2(p)e^{2t}>0,\,\forall\, t
\eee 
which by the completeness of \be M\ee implies that the integral curve \be\gamma\ee is regular and well-defined for all 
\be t\in \real{}\tchar{.}\ee Along the same curve \be\gamma\ee we now define a new function
\bee
\tilde f(t):=(|A|^2-|P*A|^2)(\gamma(t))
\eee
and from (\ref{extra 8}) we derive the evolution equation
\bee
\dt\tilde f=-2\tilde f
\eee
and
\bee
(|A|^2-|P*A|^2)(\gamma(t))=(|A|^2-|P*A|^2)(p)e^{-2t}\tchar{.}
\eee
This is the typical behaviour of a cone.
If 
\bee(|A|^2-|P*A|^2)(p)\ne 0\tchar{,}\eee 
then \be(|A|^2-|P*A|^2)(\gamma(t))\ee becomes unbounded as \be t\to-\infty\tchar{.}\ee
This contradicts the boundedness of the second fundamental form. So \be|A|^2=|P*A|^2\ee and (\ref{extra 6}) is valid.
The multiplicities of the two eigenvalues of \be P\ee are constant on \be M\ee since \be\nabla\dd Pij=0\tchar{.}\ee 
Let \be r\ee denote the multiplicity of \be\lambda=1\tchar{,}\ee i.e. \be r\ee is the rank of the second fundamental form
\be A^\nu=P/|H|\tchar{.}\ee We define the two distributions
\bcoleqn
\sE_pM&:=&\{V\in T_pM:P_i^jV^i=V^j\},\nonumber\\
\sF_pM&:=&\{V\in T_pM:P_i^jV^i=0\}\tchar{.}\nonumber
\ecoleqn
so that \be T_pM=\sE_pM\oplus \sF_pM\tchar{.}\ee 
In addition, for the nullspace \be \sE_pM\ee we have
\be\sE_pM\subset\text{ker}(\theta)\ee
since (\ref{extra 9}) implies
\bcoleq\label{extra 2}
\theta(V)=\theta_jV^j=\theta_jP_i^jV^i=0,\,\forall\, V\in \sE_pM\tchar{.}
\ecoleq
Let \be e_{p}\in\sE_ pM, f_p\in \sF_pM\ee and let \be\gamma(t)\ee be a curve in \be M\ee with \be\gamma(0)=p\tchar{.}\ee If 
\be e(t)\ee with \be e(0)=e_p\ee and \be f(t)\ee with \be f(0)=f_p\ee are varying along \be \gamma\ee by parallel transport 
we obtain from \be\nabla P=0\ee
\bee
\dt|Pe-e|^2=0=\dt|Pf|^2\tchar{,}
\eee
so that \be e(t)\in\sE_{\gamma(t)}M\ee and \be f(t)\in\sF_{\gamma(t)}M\tchar{.}\ee In particular for vector fields
\be f_1,f_2\in\sF M\ee and \be e_1,e_2\in \sE M\ee we have 
\bcoleqn
\nabla_{e_1}e_2,\,\,\,[e_1,e_2]\in\sE M,\label{extra 5a}\\
\nabla_{f_1}f_2,\,\,\,[f_1,f_2]\in\sF M\tchar{,}\label{extra 5b}
\ecoleqn
so that \be \sE M, \sF M\ee are invariant under parallel transport and in particular involutive.
This means that for each \be p\in M\ee there exist two integral leaves \be\sE_p, \sF_p\ee 
intersecting orthogonally in \be p\tchar{.}\ee 
By Lemma \ref{lemma 1}, we have for any normal vector
\be V\ee perpendicular to \be \nu\ee
\bee
\langle P, A^V\rangle=0\tchar{,}
\eee
where
\bee
A^V:=\langle V,A\rangle
\eee
is the second fundamental form w.r.t. \be V\tchar{.}\ee
At some point \be p\in M\ee let us choose an orthonormal basis \be e_1,\dots, e_r\ee of \be\sE_pM\tchar{.}\ee
Then 
\bcoleqn\label{extra 3}
&&\text{tr}_{\sE}(A^V):=\sum_{i=1}^r\langle V,A(e_i,e_i)\rangle=\langle P,A^V\rangle=0,\nonumber\\
&&\forall\, V\in N_pM\,\tchar{\text{ with }}\,\langle H,V\rangle=0\tchar{.}
\ecoleqn
{F}rom (\ref{extra 2}) we conclude that each integral leaf of \be \sE M\ee is an \be r\ee{- dimensional} submanifold of the sphere
\be S^{n-1}(\rho)\subset\real{n}\ee of radius \be\rho=|F|\tchar{,}\ee the radius depending on the leaf (see figure \ref{fig 1}).
We claim that the leaves \be \sF_p\ee
are \be(m-r)\ee{- dimensional} affine subspaces of \be\real{n}\ee and that these affine subspaces are 
parallel for any \be p,p'\ee in the same connected component of \be M\tchar{.}\ee
To see this, let \be q\in \sF_p\ee be arbitrary and note that \be \dd Aij=P_i^k\dd Akj\ee implies
\bcoleq\label{vanish}
A(f,X)=0,\,\forall\, f\in\sF_qM=T_q\sF_p,\,X\in T_qM\tchar{.}
\ecoleq
The normal space of \be \sF_p\ee at \be q\ee decomposes into 
\bee N_q\sF_p=N_qM\oplus\sE_qM\tchar{.}\eee 
By (\ref{vanish}) the second fundamental form of \be\sF_p\ee vanishes w.r.t. any normal vector belonging to \be N_qM\ee
and by (\ref{extra 5b}) also for all normal vectors belonging to \be\sE_qM\tchar{.}\ee 
This shows that the integral leaves of \be \sF M\ee are affine subspaces of \be\real{n}\ee
of dimension \be m-r\tchar{.}\ee Now we fix a point \be p\in M\ee and suppose that \be p'\in\sE_p\ee is arbitrary. Let
\be\gamma\subset\sE_p\ee be any smooth curve with \be\gamma(0)= p\ee and \be\gamma(1)= p'\tchar{.}\ee Suppose \be f_p\in\sF_p\ee is some
vector and \be f(t)\ee is varying along \be\gamma\ee by parallel transport with \be f(0)=f_p\tchar{.}\ee 
Considering \be f(t)\ee as a vector field along
\be\gamma\ee in \be\real{n}\ee we compute \be\dt f=\nabla_{\dot c}f+(\dt f)^\perp=(\dt f)^\perp=A(\dot c,f)\tchar{.}\ee
But then by (\ref{vanish}) \be\dt f=0\tchar{.}\ee Since \be\sF M\ee is invariant under parallel transport we conclude
that \be f_p\in \sF_{p'}M\ee as well and since \be f\ee was arbitrary we must have \be\sF_pM=\sF_{p'}M\tchar{.}\ee
To see that this relation holds for all \be p, p'\ee in the same connected component of \be M\ee we proceed as follows. The leaf
\be\sE_p\ee is contained in a linear subspace of dimension \be n-m+r\ee which is orthogonal to \be\sF_pM\tchar{.}\ee
Then we define an embedding
\bcoleqn
\tilde F:\sE_p\times\real{m-r}\to\real{n},\nonumber\\
\tilde F(p',X)=(p',X).\nonumber
\ecoleqn
Since the integral leaves of \be \sF M\ee through \be p'\in\sE_p\ee coincide with the affine subspace w.r.t. \be \sF_{p}M\ee through
\be p'\tchar{,}\ee the \be m\ee{- dimensional} image of \be\sE_p\times\real{m-r}\ee under \be\tilde F\ee is contained in
the \be m\ee{- dimensional} submanifold \be M\ee and by the completeness of both \be\sE_p\ee and \be M\ee we
conclude that \be\tilde F(\sE_p\times\real{m-r})\ee is a connected component of \be M\tchar{.}\ee
Since the integral leaves of \be\sE M\ee are contained in spheres of radii \be\rho=|F|\ee and 
\bee|F|^2=|F^\perp|^2+|\theta|^2=|H|^2+|\theta|^2
\eee
we observe that \be|\theta|^2\ee must be constant along the leaves of \be\sE M\ee since \be H\ee is parallel and hence \be|H|\ee
is constant. Moreover, from (\ref{extra 2}) we conclude that the projection \be \pi:M\to \sE_p\ee maps the points
\be q\ee in the integral leaves of \be\sE M\ee with \be|q|=|F|\ee to points \be\tilde q\in\sE_p\ee with 
\be|\tilde q|_{\sE_p}=|F^\perp|=|H|\ee and since \be H=-F^\perp\ee we see that \be \tilde M:=\sE_p\ee must be 
a minimal submanifold of the sphere
\be S^{n-m-1+r}{(|H|)}\tchar{.}\ee On the other hand \be|H|^2=\uu gij\dd Pij=r\ee since all eigenvalues of \be P\ee are
either \be 0\ee or \be 1\ee and \be r\ee is the multiplicity of the latter. Thus \be|H|=\sqrt{r}\ee and 
each connected component of \be M\ee is isometric to \be \tilde M\times\real{m-r}\ee with a minimal submanifold
\be \tilde M\subset S^{n-m-1+r}(\sqrt{r})\tchar{.}\ee \\
\begin{figure}[htbp]
    \includegraphics[width=.5\hsize]
      {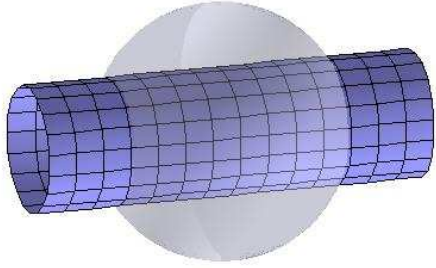}
    \caption{The sphere intersects the cylinder \be M\ee in two leaves of the distribution \be\sE M\ee which are the circles
on the cylinder. The straight lines perpendicular to the circles are the leaves belonging to \be\sF M\tchar{,}\ee}
\label{fig 1}
\end{figure}

\item[(ii)]The remaining case is \be\nabla^\perp H(p_0)\neq 0\ee at a point \be p_0\in M\tchar{.}\ee Then, as in the compact
case, we find \be|P|^2=|H|^4\ee globally. In addition, \be P\ee admits only one non-zero eigenvalue \be\lambda=|H|^2\ee  and 
\be\frac{\nabla|H|}{|\nabla|H||}\ee spans the eigenspace. Then Lemma \ref{lemma 1} implies the relations
\bcoleq\label{ex -1}
P_i^k\dd Akj=|H|\dd Pij\,\nu,\quad P_i^j\,\nabla^i|H|=|H|^2\nabla^j|H|\tchar{.}
\ecoleq
Let us define
\bee
\cM:=\{p\in M:\nabla H(p)\neq 0\}\tchar{.}
\eee
{F}ollowing Huisken's \cite{huisken2} argument closely, we choose a connected component
\be U\subset\cM\ee and consider the distributions 
\bcoleqn
\sE_pU&:=&\{X\in T_pU:p\in U,\,PX=|H|^2X\}\tchar{,}\nonumber\\
\sF_pU&:=&\{X\in T_pU:p\in U,\,PX=0\}\tchar{.}\nonumber
\ecoleqn
In contrast to part (i), these distributions are defined a-priori only on \be\cM\tchar{.}\ee
As in the first part, they  are invariant under parallel transport and again (\ref{extra 5a}), (\ref{extra 5b})
are valid for them so that we may apply Frobenius' theorem. We claim that the integral leaves of \be\sF U\ee
are affine subspaces of dimension \be m-1\tchar{.}\ee This is equivalent to showing that
\bee
\dAcirc{ij}:=\dd Aij-\frac{1}{|H|}\dd Pij\,\nu
\eee
vanishes completely. From (\ref{mean}) we observe
\bcoleq\label{ex 0}
|H|\nabla_i|H|=\theta^k\dd Pki
\ecoleq
and then with (\ref{ex -1})
\bee
\theta\left(\nabla|H|\right)=\frac{|\nabla|H||^2}{|H|}\,\tchar{.}
\eee
Let \be\dthetacirc{}\ee be the projection of \be\theta\ee onto \be\sF U\tchar{,}\ee i.e.
\bee
\dthetacirc{i}=\theta_i-\frac{\theta\left(\nabla|H|\right)}{|\nabla|H||^2}\,\nabla_i|H|=\theta_i-\frac{1}{|H|}\nabla_i|H|\tchar{.}
\eee
A straightforward computation shows
\bcoleqn
\uthetacirc{k}\dAcirc{ki}
&=&\left(\theta^k-\frac{1}{|H|}\nabla^k|H|\right)\left(\dd Aki-\frac{1}{|H|}\dd Pki\,\nu\right)\nonumber\\
&\overset{(\ref{ex 0})}=&\theta^k\dd Aki-\nabla_i|H|\,\nu-\frac{1}{|H|}\nabla^k|H|\dd Aki\nonumber\\
&&+\,\frac{1}{|H|^2}\nabla^k|H|\dd Pki\,\nu\nonumber\\
&\overset{(\ref{ex -1})}=&\theta^k\dd Aki-\frac{1}{|H|^3}\nabla^j|H|P_j^k\dd Aki\nonumber\\
&\overset{(\ref{ex -1})}=&\theta^k\dd Aki-\frac{1}{|H|^2}\nabla^j|H|\dd Pij\nu\nonumber\\
&\overset{(\ref{ex -1})}=&\theta^k\dd Aki-\nabla_i|H|\nu\nonumber\\
&\overset{(\ref{mean})}=&0\tchar{.}\label{ex 5}
\ecoleqn
Then
\bcoleqn
0&=&\nabla_l^\perp\left(\uthetacirc{k}\dAcirc{ki}\right)\nonumber\\
&\overset{(\ref{ss 2})}=&\dAcirc{li}-P_l^k\dAcirc{ki}+\frac{1}{|H|^2}\nabla_l|H|\nabla^k|H|\dAcirc{ki}\nonumber\\
&&-\,\frac{1}{|H|}\nabla_l\nabla^k|H|\dAcirc{ki}+\uthetacirc{k}\nabla_l^\perp\dAcirc{ki}\nonumber\\
&\overset{(\ref{ex -1})}=&\dAcirc{li}-\frac{1}{|H|}\nabla_l\nabla^k|H|\dAcirc{ki}+\uthetacirc{k}\nabla_l^\perp\dAcirc{ki}\tchar{.}
\nonumber
\ecoleqn
Since \be\sF U\ee is invariant under parallel transport, we have
\bee
\nabla_{f_1}\nabla_{f_2}|H|=f_1\left(\langle\nabla|H|,f_2\rangle\right)-\langle\nabla|H|,\nabla_{f_1}f_2\rangle=0
\eee
for all \be f_1,f_2\in\sF U\tchar{.}\ee Therefore we may write
\bee
\nabla_i\nabla_j|H|=\frac{\Delta|H|}{|H|^2}\dd Pij
\eee
and since (\ref{ex -1}) implies \be P_i^k\dAcirc{\hspace{-1pt}kj}=0\tchar{,}\ee
we get from the formula above
\bcoleq
\uthetacirc{k}\nabla_l^\perp\dAcirc{ki}=-\dAcirc{li}\tchar{,}\label{ex 2}
\ecoleq
in particular both sides in (\ref{ex 2}) are symmetric in \be i\ee and \be j\tchar{.}\ee
{F}rom (\ref{ss 2}), \be\nabla|H|(p)\in\sE_p U\tchar{,}\ee \be\dthetacirc{}(p)\in\sF^*_pU\ee and since
\be\sF U\ee is invariant under parallel transport we derive
\bcoleq\label{ex 6}
\eta\left(\nabla|\dthetacirc{}|^2\right)=2\Bigl\langle\eta,\dthetacirc{}\Bigr\rangle,\quad\forall\, \eta\in\sF^*_p U\tchar{.}
\ecoleq
We may now proceed similarly as in part (i). If \be\dthetacirc{}(p)=0\tchar{,}\ee then (\ref{ex 2}) implies \be\dAcirc{}=0\tchar{.}\ee
If \be\dthetacirc{}(p)\neq 0\tchar{,}\ee we choose \be\eta=\dthetacirc{}\ee in (\ref{ex 6}) and consider the integral curve \be\gamma\ee
w.r.t. \be\dthetacirc{}\ee starting at \be p\ee which as before proves that the integral curve is regular and contained in \be U\ee
for all \be t\in\real{}\tchar{.}\ee Applying formula (\ref{ex 2}) to \be\gamma\ee we get exponential growth in \be t\ee for the quantity
\be|\dAcirc{}|^2\ee if \be|\dAcirc{}(p)|\neq 0\ee which by the boundedness of the second fundamental form would be a contradiction.
Therefore \be\dAcirc{}=0\ee everywhere and as in the first part we obtain for some fixed point \be p\in U\ee
that the integral leaf \be \sE_p U\ee is part of an Abresch-Langer curve \be\Gamma\ee that extends smoothly up to the boundary
of \be U\ee and that all integral leaves \be\sF_q U\ee for
\be q\in U\ee are parallel affine subspaces. So as above any connected component \be U\subset\cM\ee 
is a product of a part of an Abresch-Langer curve and some \be R^{m-1}\tchar{.}\ee To extend this to all of \be M\ee we have to consider
the points, where \be\nabla|H|=0\tchar{.}\ee We claim that the critical points of \be|H|\ee are transversally isolated along 
\be\sE U\tchar{.}\ee By this we mean the following. The integral leaves \be\sE_pU\ee of \be U\ee are geodesics in \be M\ee and
by the completeness can be extended arbitrarily. So we obtain a complete geodesic \be\gamma\subset M\ee that coincides in \be U\ee
with a portion of an Abresch-Langer curve \be\Gamma=\sE_pU\ee and which intersects the boundary of \be U\ee orthogonally.
Suppose the fraction of  \be\gamma\subset U\ee is parametrized by \be\gamma:(0,t_0)\to U\ee and that 
\be\lim_{t\to t_0}\nabla |H|(\gamma(t))=0\tchar{.}\ee Then (\ref{ss 6}) implies 
\bcoleq\label{mm}
\lim_{t\to t_0}\Delta|H|=|H|-|H|^3
\ecoleq
On the other hand \be |H|(\gamma(t))\ee 
coincides in \be U\ee with the curvature \be k\ee of the Abresch-Langer curve
\be\Gamma\ee and \be\lim_{t\to t_0}k(t)=k(t_0)=k_0\ee is a critical value
of the curvature of \be\Gamma\tchar{.}\ee Since \be k_0=1\ee is the critical value of the round circle, lemma \ref{abl}  and (\ref{mm})
would imply that, if \be\Delta|H|(t_0)=0\tchar{,}\ee then \be\Gamma\ee is a fraction of \be S^1\ee and \be\nabla k=0=\nabla|H|\ee on all of 
\be\Gamma\tchar{.}\ee Therefore \be k_0\neq 1\ee and \be\Delta |H|(t_0)\neq 0\tchar{.}\ee  This implies that \be\nabla |H|\neq 0\ee
immediately after passing \be \gamma(t_0)\ee in direction of \be \gamma\tchar{.}\ee In particular, the boundary of each connected component
\be U\subset\cM\ee consists of at most two affine subspaces of dimension \be m-1\ee that are parallel to each of the leaves
\be\sF_qU,\,q\in U\tchar{.}\ee If \be U_1, U_2\ee are two connected components with \be\bar U_1\cap\bar U_2\neq\emptyset\tchar{,}\ee
then the integral leaves \be\sF_{q_1}U_1,\,\sF_{q_2}U_2\ee are parallel.
We may thus join two such connected components \be U_1, U_2\ee along their common boundary and obtain two portions of Abresch-Langer
curves \be\Gamma_1, \Gamma_2\ee that both intersect the common boundary of \be U_1, U_2\ee orthogonally and which meet at a boundary point,
for instance at \be\gamma(t_0)\tchar{.}\ee
If \be\Gamma_1\ee and \be\Gamma_2\ee would be portions of different Abresch-Langer curves, then again by 
lemma \ref{abl} we get a contradiction, since then
the critical values of the two curvature functions are different but both satisfy (\ref{mm}). By the completeness of \be M\ee we may
extend this process until we get \be M=\Gamma\times\real{m-1}\tchar{.}\ee
\end{itemize}
\end{pproof}

\nocite{*}
\bibliographystyle{alpha}
\bibliography{refs}

\end{document}